\documentclass[]{article}

\usepackage{amssymb} 
\usepackage{amsmath} 
\usepackage{latexsym}
\usepackage{graphicx}

\usepackage{enumerate}

\usepackage{tikz}
\usetikzlibrary{calc}
\usepackage{pgf,tikz}
\usetikzlibrary{arrows}
\usepackage{hyperref}
\usepackage{amscd}
\usepackage{relsize}

\usepackage{tikz-cd}

\usepackage{ color, amsmath, amssymb, amsfonts, amscd, graphicx,latexsym, hyperref}
\usepackage{caption, subcaption}
\usepackage[pdf]{pstricks}

\usepackage[affil-it]{authblk}

\usepackage{amsthm}

\newtheorem{theorem}{Theorem}

\newtheorem{proposition}{Proposition}

\newtheorem{remark}{Remark}

\newtheorem{lemma}{Lemma}

\newtheorem{definition}{Definition}

\newtheorem{corollary}{Corollary}

\newtheorem*{example}{Example}

\title{Complete Flat Cone Metrics on Punctured Surfaces}
\author{ \.{I}smail Sa\u{g}lam \thanks{Electronic address: \texttt{isaglamtrfr@gmail.com}  }}
\affil{Adana Science and Technology University  }

\date{}
\begin{document}

\maketitle

\begin{abstract}
We prove that each complete flat cone metric on
a surface, perhaps with boundary and punctures, can be triangulated with finitely many types of triangles. We derive Gauss-Bonnet formula for this kind of cone metrics. In addition, we prove that 
each free homotopy class of paths has a geodesic representative. 
\end{abstract}
\thanks{}
\tableofcontents
\section{Introduction}
\thanks{}
Flat cone metrics  appear in several areas of mathematics. For example, they are 
studied in Teichm\"{u}ller theory through 
quadratic differentials,  and in dynamics of billiard tables \cite{teichmuller}, \cite{dynamics}. These objects are also interesting for their own sake. Classification of certain families of these metrics may give rise to interesting results in several areas of mathematics such as hypergeometric functions, (real and complex) hyperbolic geometry \cite{delign-mostow}, \cite{thurston}, \cite{bavard-ghys}, \cite{fillastre}, \cite{fillastre-izmestiev}, \cite{rivin}. In addition, regarding 
combinatorial triangulations or quadrangulations as cone metrics as in \cite{thurston}, one can parametrize certain families of dessins d'enfants. See \cite{ayberk-muhammed},
\cite{ismail-muhammed}, \cite{ayberk}, \cite{ayberk-lattes}, \cite{allcock}.  

Flat cone metrics on compact surfaces have been studied well.  We know that there is a length minimizing path between any two point of such a surface. Also, each free homotopy class of loops on a compact surface with a flat metric contains a length minimizing geodesic.
Indeed, these properties follow from general theory of  length spaces \cite{Gromov}, \cite{dimitri}.  Furthermore, Gauss-Bonnet formula holds for these surfaces, and they can be triangulated with finitely many triangles.
See \cite{Tro-enseign}, \cite{Tro-handbook}, \cite{Tro-compact}.

Teichm\"{u}ller theory is related with the theory of cone metrics in a natural way. Let $S$ be a closed, 
orientable surface, $x_1, \dots, x_n \in S$. Pick  $a_1, \dots a_n \in \mathbb{R}$ so that $\sum_{i=1}^{n}a_i=2\pi\chi(S)$, where $\chi(S)$
is Euler characteristics of $S$. Consider the curvature divisor 

$$D=\sum_{i=1}^n a_ix_i.$$
\noindent It is known that  each conformal class on $S$ includes a flat metric with $n$ singular points of divisor $D$. Furthermore, this metric is unique 
up to homothety. See the papers of Troyanov \cite{Tro-enseign}, \cite{Tro-compact} in case  $a_i<2\pi$, for each $1 \leq i \leq n$. Otherwise, see the paper of Hulin and Troyanov \cite{Tro-open}.

Flat surfaces \textit{regular} punctures have been also studied well. By a regular puncture on a flat surface, we mean a puncture which has a neighborhood isometric to that of point at infinity of a cone. Gauss-Bonnet formula holds for the surfaces with regular punctures. Also, there is length minimizing geodesic  in any homotopy class of loops in such a surface. In addition, these surfaces may be triangulated with finitely many types of triangles.

Our objective  is to verify 
that any complete flat metric on a given surface, with regular or irregular punctures, has  above mentioned properties. Let $\bar{S}$ be a surface with a complete flat cone metric. We   summarize the results of the present paper.
\begin{enumerate}
\item 
In Section \ref{triangulation-subsection}, we show that $\bar{S}$ can be triangulated with finitely many types of triangles. 
\item
In Section \ref{Gauss-Bonnet}, we show that a variant of Gauss-Bonnet formula holds for $\bar{S}$. 
\item
In Section \ref{free-homotopy}, we show that each loop on $\bar{S}$ has a geodesic representative in its free homotopy class.
	\end{enumerate}

We want to study complete flat metrics at the highest level of generality. The surfaces that we consider  are of finite type, may have punctures and boundary.
We do not omit the surfaces having punctured boundary components from our discussion. Therefore, we start with introducing a convenient  notation.

\subsection{Notation}
\label{notation}

\begin{definition}
	Let $S$ be a compact, connected topological  surface perhaps with boundary $B$.
	Let $\mathfrak{l}, \mathfrak{p},\mathfrak{l'},\mathfrak{p'}$ be finite disjoint subsets of $S$ so that 
	\begin{itemize}
		\item
		$\mathfrak{l}$ and $\mathfrak{p}$ are  subsets of the \textit{interior} of $S$,
		\item
		$\mathfrak{p'}, \mathfrak{l'}$ are  subsets of $B$.
		
	\end{itemize}

	An element in $\mathfrak{l}$ will be called  \textit{labeled interior}  point.  An element in $\mathfrak{p}
	$ will be called punctured interior point. Other points in interior of $S$ called \textit{ordinary interior points} An element in $B$ will be called  boundary point. An element
	in $\mathfrak{l'}$ will be called a \textit{labeled boundary} point. An element in $\mathfrak{p'}$
	will be called  \textit{punctured boundary} point.
	Other points in boundary will be called \textit{ordinary boundary points}. 
	A \textit{doubly labeled surface}, shortly DL surface, is the tuple 
	$$ (S,B, \mathfrak{l},\mathfrak{p}, \mathfrak{l'}, \mathfrak{p'} )$$
\end{definition}

Also we will use the following notation:
\begin{itemize}
	\item
	$S_B=S-B$.
	\item
	$S_{\frak{l}}=S-\frak{l}$
	\item
	$S_{B,\frak{l}}=S-(B\cup \frak{l})$
	\item
	\dots
\end{itemize}

We will denote a doubly labeled surface $ (S,B, \mathfrak{l},\mathfrak{p},\mathfrak{l'}, \mathfrak{p'} ) $ as $S^L$. Underlying compact surface of $S^L$ will simply be denoted as $S$. 

DL surfaces can be considered as punctured surfaces, with puncture set  $\frak{p} \cup \mathfrak{p'}$. Indeed, $S_{\frak{p,p'}}$ is the punctured surface that we consider. Observe that 
punctured and labeled points may lie in boundary also.

\section{Flat DL  surfaces}
Flat compact surfaces can be triangulated with finitely many triangles. For non-compact surfaces, we need to modify the definition of triangulation. The reason for this is that punctured surfaces may require infinitely many triangles and arbitrary triangulations possibly  induce  
non-complete cone metrics. 
\begin{definition}
	\label{triangulation}
	A Euclidean triangulation of a DL surface $S^L$ is a set of pairs $\mathfrak{T}=\{(T_{\alpha}, f_{\alpha})_{\alpha \in A}\}$ where each $T_{\alpha}$ is a compact subset of
	$S_{\frak{p}, \mathfrak{p'}}$ and $f_{\alpha}:T_{\alpha} \rightarrow \mathbb{R}^2$ is a homeomorphism onto a non-degenarete triangle $f_{\alpha}(T_{\alpha})$ in the Euclidean plane. $T_{\alpha}$ is called a triangle. Let $e$ be a subset of $T_{\alpha}$. $e$ is called an edge if 
	$f_{\alpha}(e)$ is an edge for the Euclidean triangle $f_{\alpha}(T_{\alpha})$.
	Similary, $v\in T_{\alpha}$ is called a vertex if $f_{\alpha}(v)$ is a vertex 
	of the triangle $f_{\alpha}(T_{\alpha})$. The Euclidean triangulation also satisfies the following properties:
	\begin{enumerate}
		\item
		$S_{\mathfrak{p,p'}}=\cup_{\alpha \in A}T_{\alpha}$
		\item
		If $\alpha \neq \beta$, then $T_{\alpha}\cap T_{\beta}$ is either empty or an \textit{edge}, or a \textit{vertex}.
		\item
		If $T_{\alpha}\cap T_{\beta}$ is not empty, then there is a $g_{\alpha \beta} \in \mathfrak{E}(2)$ (the group of isometries of Euclidean plane) so that $f_{\alpha}=g_{\alpha \beta}f_{\beta}$ on the intersection.
		\item
		(Local finiteness) Each compact subset of $S_{\mathfrak{p,p'}}$ intersect  with finitely many \textit{triangles}, \textit{edges} and \textit{vertices}.
		\item
		Set of triangles $f_{\alpha}(T_{\alpha})$ consisits of finitely many isometry 
		classes of Euclidean triangles. 
	\end{enumerate}
\end{definition}
Observe that our definition is a generalization of the one given in \cite{Tro-handbook}. We just added two more conditions: $\textup{(4),(5)}$. Note that 
Euclidean triangulations on compact surfaces always have these properties.
We will show that these triangulations induce complete flat cone metrics on
DL surfaces and DL surfaces with complete cone metrics can be triangulated. See Proposition \ref{proposition-triangle} and Theorem \ref{theorem-triangle}.

\paragraph{The notions of angle and curvature for DL surfaces having Euclidean triangulations}

\begin{definition}
	Let $S^L$ be a DL surface together with a Euclidean triangulation. Let  $x$ be a vertex in $S_B$. $x$ is called a point having angle  $\theta$ if
	\begin{align}
	\theta=\theta(x)=\sum_{j=1}^k \phi_{j},
	\end{align}
	where $\phi_1,\dots, \phi_k$ are angles of the triangles incident to $x$, at the vertex $x$.  The curvature at $x$ is
	\begin{align}
	\kappa=\kappa(x)=2\pi - \theta(x).
	\end{align}
	Similarly, let  $y$ be a vertex in $B- \mathfrak{p'}$. $y$ is called  a point having angle  $\theta$ if
	\begin{align}
	\theta=\theta(y)=\sum_{j=1}^{r} \phi'_j,
	\end{align}
	where $\phi'_1,\dots, \phi'_r$ are angles of the triangles incident to $y$, at the vertex $y$. The curvature at $y$ is
	\begin{align}
	\kappa=\kappa(y)=\pi - \theta(y).
	\end{align}
	Curvature at the points which are not vertices is defined to be $0$. If a point, either on the boundary
	or not, has curvature $0$, then it is called non-singular. Otherwise it is called singular.
\end{definition}

\begin{definition}
	A flat doubly labeled (FDL) surface  $(S^L,\mathfrak{T})$ is a DL surface $S^L$ together with an Euclidean triangulation $\mathfrak{T}$ such that its set of singular points is  $\mathfrak{l}\cup \mathfrak{l'}$. 
\end{definition}
In Section \ref{Gauss-Bonnet}, we will extend the notions of   the curvature and  the angle to the punctured interior and punctured boundary points.

\subsection{Induced length structure }

An FDL surface $(S^L,\frak{T)}$ has natural area measure which coincides with the  2 dimensional Lebesque measure at each triangle $T_{\alpha}$. Also, 
as in \cite{Tro-handbook}, we can define the length $l(c)$ of a curve 
$c: \ [a,b] \rightarrow S_{\frak{p,p'}} $ ($a,b \in \mathbb{R}, \ a<b$) as follows:

\begin{itemize}
	\item 
	If $c$ is contained in a triangle $T_{\alpha}$ of $\frak{T}$, then 
	$l(c)$ is its Euclidean length.
	\item
	If $c$ is concatenation of two curves $c_1$ and $c_2$, then $l(c)=l(c_1)+l(c_2)$.
\end{itemize}

When there is no risk of confusion we will refer curves on $S_{p,p'}$ as
curves on $S^L$. Also we will use the notation 
$[a,b] \rightarrow S^L$ instead of $[a,b] \rightarrow S_{p,p'}$.

\begin{lemma}
	\label{admissible}
	Let $S^L$ be a FDL surface. Any two points in $S_{\mathfrak{p,p'}}$ can be joined by a
	curve of finite length.
	
	\begin{proof}
		Take two points $x,y \in S_{\frak{p,p'}}$ and a curve $c: \ [0,1] \rightarrow  S^L$ joining them. Since image of the curve is compact it is 
		contained in a finite number of triangles. One can easily construct 
		a finite length curve joining $x$ and $y$ which lies in the union of these triangles.
	\end{proof}
	
\end{lemma}

Consider the following function $d:\ S_{\frak{p,p'}}\times S_{\frak{p,p'}}\rightarrow \mathbb{R}$
\begin{align}
d(x,y)=\inf \{l(\alpha): \textup{ $\alpha$  is  a curve joining  x to y } \}.
\end{align}

\begin{proposition}\label{metric}
	$d$ is a metric on $S_{\mathfrak{p,p'}}$:
	\begin{enumerate}
		\item
		$d(x,x)=0$ for all $x \in S_{{\frak{p,p'}}}$.
		\item
		$\infty >d(x,y)>0$, when $x\neq y$.
		\item
		$d(x,y)=d(y,x)$ for all $x,y \in S_{\mathfrak{p,p'}}$.
		\item
		$d(x,y)+d(y,z)\geq d(x,z)$ for all $x,y,z \in S_{\mathfrak{p,p'}}$.
	\end{enumerate}

	\begin{proof}
		$\textup{(1), (3) and (4)}$ are obvious. $\textup{(2)}$ follows from 
		local finiteness property of the Euclidean triangulations and Lemma \ref{admissible}.
	\end{proof}
\end{proposition}

If there is no risk of confusion, we will refer this metric as a metric 
on $S^L$ instead of a metric on $S_{\mathfrak{p,p'}}$.

\begin{figure}
	\begin{center}
		\includegraphics[scale=0.45]{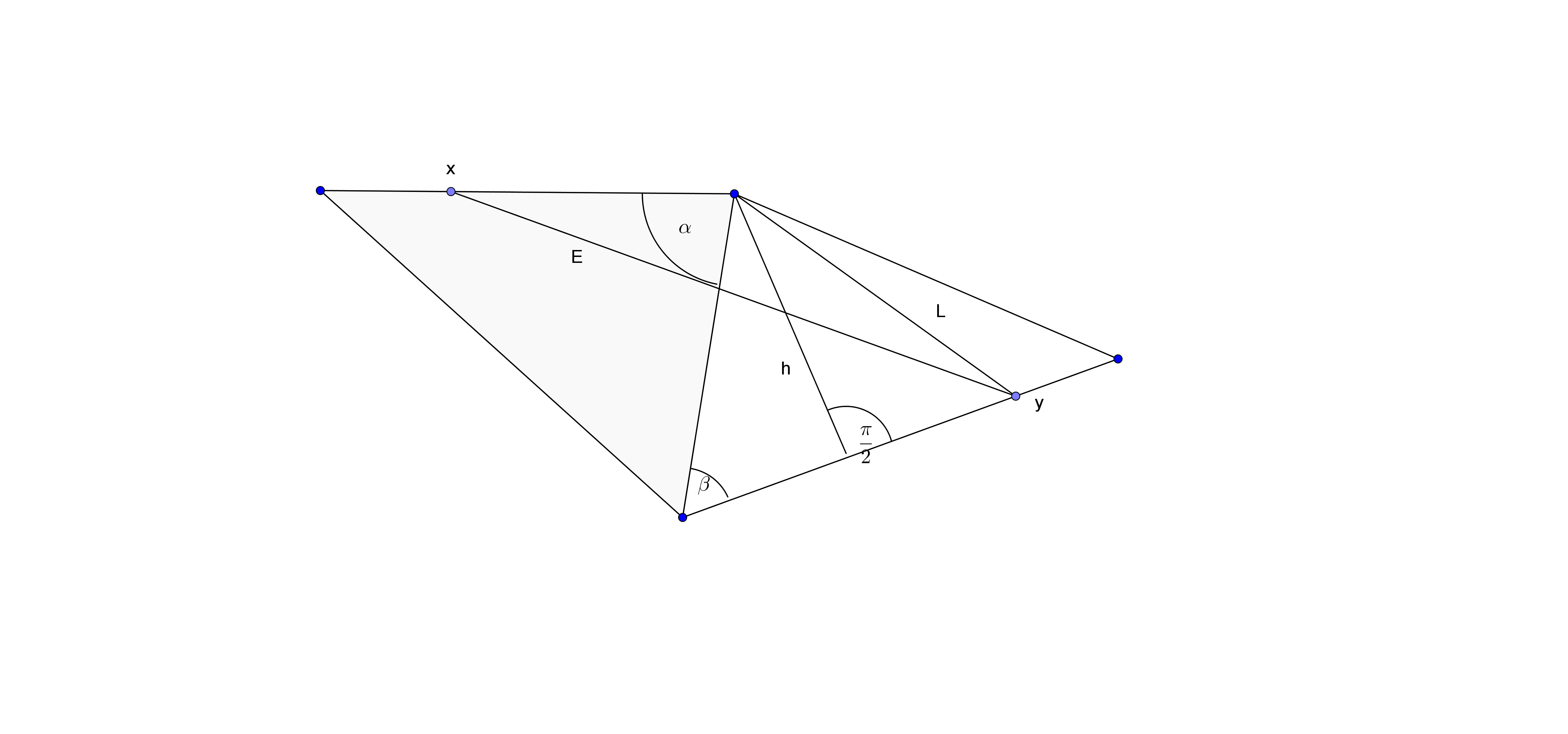}
	\end{center}
	\caption{Length of $E$ is greater than or equal to length of $h$.}
	\label{hinge}
\end{figure}

Now we state an elementary lemma from Euclidean geometry. See Figure \ref{hinge}.
\begin{lemma}
	\label{hinge-lemma}
	Let $H$ be hinge of two triangles. Let $E$ be a line segment on $H$
	joining two edges which are not adjacent. The length of $E$ is greater than or equal to one of the altitudes of the triangles.
	\begin{proof}
We assume that $\alpha\geq\beta$. It follows
	that the length of $E$ is greater than or equal to thelength of $L$, which 
	is greater than or equal to the length of $h$.
	\end{proof}
\end{lemma}

\begin{figure}
	\begin{center}
		\includegraphics[scale=0.45]{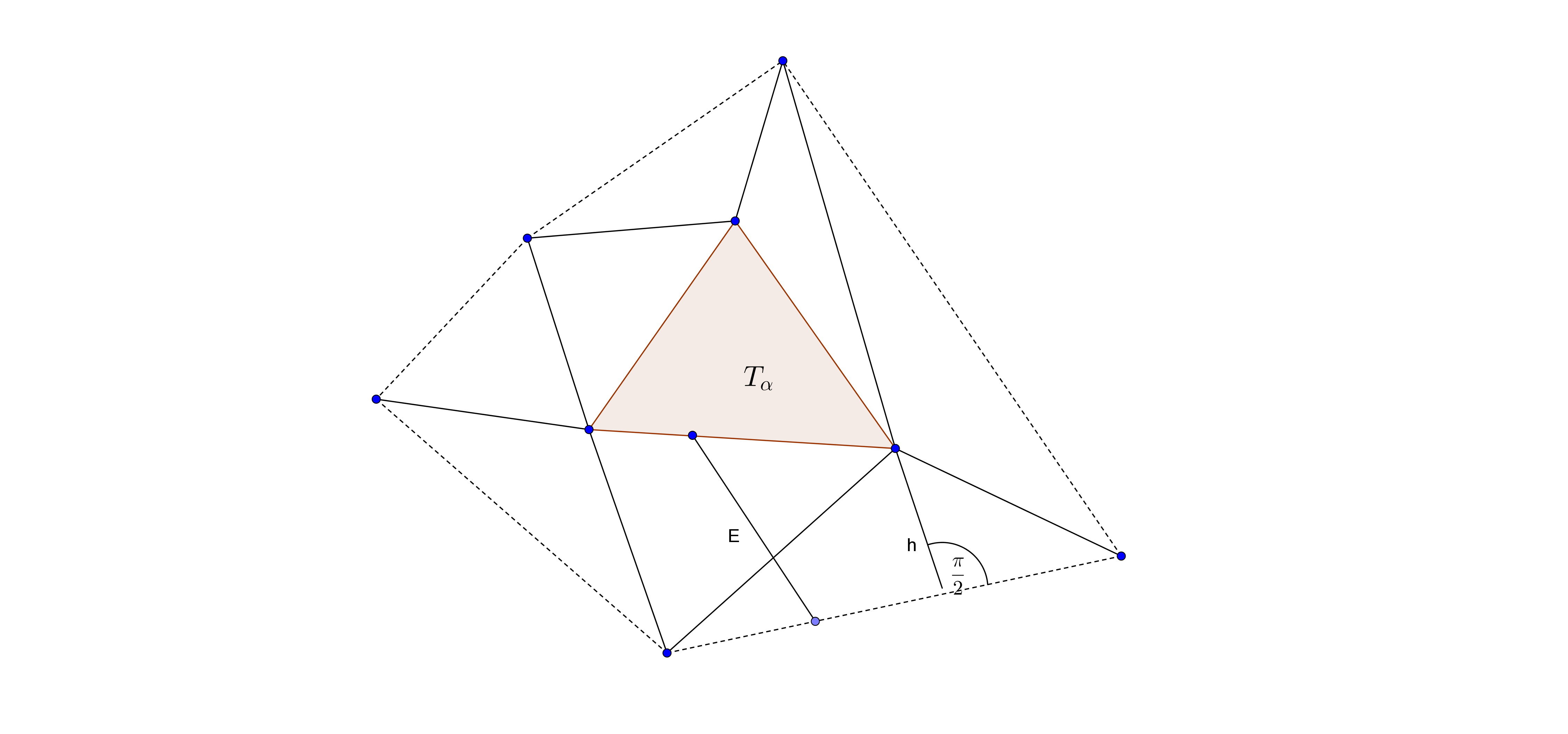}
	\end{center}
	\caption{The triangular neighborhood $U$ of $T_{\alpha}$. Dashed line segments correspond to $V$.}
	\label{complete-figure}
\end{figure}

\begin{lemma}
	\label{plane-fact}
	Let $(S^L,d)$ be an FDL surface with induced metric $d$. Let $T_{\alpha}$ be a triangle on it, and $\delta$ be  the minimum of altitudes of the triangles on it. Let $x\in T_{\alpha}$, and $y$ be a point which is not in $T_{\alpha}$ or triangles intersecting $T_{\alpha}$.
	It follows that $$d(x,y)>\delta.$$
	\begin{proof}Let
$$U=\{x\in S_{\mathfrak{p,p'}}: \ x \ \textup{is in one of the triangles intersecting with $T_{\alpha}$}\}.$$

Consider the following subset of $U$:

$$V=\{x\in U: x \ \textup{isin  an edge  which does not intersect with}\ T_{\alpha} \}.$$
 See Figure \ref{complete-figure}. Note that 
 
 \begin{enumerate}
 	\item 
 	if $S^L\neq U$, then $S^L-V$ is disconnected,
 	\item 
 	if $T$ is a triangle of $U$ so that its edge $e$ is in $V$, then distance between a point in $T_{\alpha}$ and a point in $e$ is greater than or equal to $\delta$.  See  Lemma \ref{hinge-lemma}.
 	
 \end{enumerate}
Take a curve joining $y $ to $x$. It follows that the curve and $V$
intersect. Hence $d(x,y)$ is strictly greater than distance between the sets $V$ and $T_{\alpha}$. Since the distance between $V$
and $T_{\alpha}$ is greater than or equal to $\delta$, it follows that $d(x,y)>\delta$.
	\end{proof}
\end{lemma}

\begin{proposition}
	\label{proposition-triangle}
	$(S^L,d)$ is complete metric space.
	\begin{proof}
		Let $x_1,\dots, x_n,\dots$ be a Cauchy sequence in $S_{\frak{p},\mathfrak{p'}}$. There exists $m \in \mathbb{Z}^+$ such that for all $n \geq m$  $d(x_n,x_m) < \delta$, where $\delta$ is minimum of the lengths of the altitudes of all triangles $T \in \frak{T}$. Let $T_{\alpha}$ be one of the triangles which 
		contains $x_m$. By locally finiteness, the set of all triangles incident to
		$T_{\alpha}$, either from a vertex or from an edge, is finite. Consider 
		the following compact set:
		$$U=\{x\in S_{\mathfrak{p,p'}}: \ x \ \textup{is in one of the triangles incident to $T_{\alpha}$}\}.$$
		If $U=S^L$, then $S^L$ is compact and so it is complete. If $U\neq S^L$, Lemma \ref{plane-fact}  implies that $U$ contains a ball of radius 
		$\delta$ around $x_m$. Hence $x_m,x_{m+1}, \dots$ are contained in $U$. Since $U$ is compact, this sequence  converges to some element in $U$. This means that 
		the sequence $(x_n)_{n=1}^{\infty}$ is convergent.
	
	\end{proof}
\end{proposition}

\begin{remark}
	$(S^L,d)$   is a length space. See \cite{Gromov}.
\end{remark}

\begin{proposition}
	\begin{enumerate} 
		\item 
		Given any two points on $S_{\frak{p},\mathfrak{p'}}$, there exists a path joining them which has minimum length.
		\item
		If $\frak{p}\cup \mathfrak{p'}$ is not empty then $\mathfrak{T}$ contains infinitely many triangles.
		\item
		If $\frak{p}\cup\mathfrak{p'}$ is not empty then $d$ is unbounded.
	\end{enumerate}
	\begin{proof}
		\begin{enumerate}
			
			\item
			This follows from Hopf-Rinow theorem for complete length spaces. 
			\item
			If there are finitely many triangles then $S_{\mathfrak{p,p'}}$ is compact.  This is impossible.
			\item
			A complete metric space together with a bounded metric is compact.
		\end{enumerate}
	\end{proof}
\end{proposition}

\subsection{Cones}

A cone having angle $\theta$, or equivalently  curvature $\kappa =2\pi-\theta$, is the set
\begin{align}
\{(r, \psi): r \in \mathbb{R}^{\geq 0}, \psi \in \mathbb{R}/\theta \mathbb{Z}\}
\end{align}
with the metric
\begin{align}
\mu=dr^2+r^2d\psi^2.
\end{align}

\noindent See \cite{Tro-enseign} for more information about cones. A cone can be considered as a  FDL   sphere with one punctured and one labeled point. The point $(0,0)$ is called vertex of the cone. 
We will denote origin of a cone by $v_{0}$. Since a cone can be considered as a piecewise flat surface, it makes sense to talk about the punctured point  or the \textit{point at infinity}.
We will denote this point as $v_{\infty}$.

\begin{definition}
	Consider a cone with angle $\theta >0$.
	\begin{itemize}
		\item
		$\kappa{(v_\infty)}=2\pi+\theta$ is called curvature at $v_{\infty}$.
		\item
		$\theta{(v_\infty)}=-\theta$ is called the angle at $v_{\infty}$.
		
	\end{itemize}
\end{definition}

\begin{remark}
	Observe that
	$\kappa(v_0)+\kappa{(v_\infty)}=4\pi$: Gauss-Bonnet formula
	for the sphere  holds.
	
\end{remark}

\noindent A cone with angle $\theta$ will be denoted by $C_{\theta}$.

\begin{figure}
	\definecolor{qqqqff}{rgb}{0,0,1}
	\definecolor{cqcqcq}{rgb}{0.75,0.75,0.75}
	\begin{tikzpicture}[line cap=round,line join=round,>=triangle 45,x=1.0cm,y=1.0cm]
	\fill[dash pattern=on 1pt off 2pt on 5pt off 4pt,color=cqcqcq,fill=cqcqcq,fill opacity=0.35] (-3.42,0.04) -- (10.72,0) -- (10.76,-2.12) -- (-3.42,-2.26) -- cycle;
	\fill[color=cqcqcq,fill=cqcqcq,fill opacity=0.25] (-1.56,2.75) -- (8.53,1.06) -- (8.51,5.22) -- cycle;
	\fill[color=cqcqcq,fill=cqcqcq,fill opacity=0.25] (8.51,5.25) -- (8.49,1.05) -- (11.26,0.59) -- (11.22,5.92) -- cycle;
	\draw [line width=2pt] (-1.56,2.75)-- (11.22,5.92);
	\draw [line width=2pt] (-1.56,2.75)-- (11.31,0.58);
	\draw [line width=2pt,domain=-3.990532339920091:10.788724508629592] plot(\x,{(-0.95--0.02*\x)/-14.66});
	\draw [line width=2pt,domain=-3.8544076057887127:11.352669835745303] plot(\x,{(-34.67--0.14*\x)/15.03});
	\draw [shift={(-1.56,2.75)}] plot[domain=-0.17:0.24,variable=\t]({1*2.16*cos(\t r)+0*2.16*sin(\t r)},{0*2.16*cos(\t r)+1*2.16*sin(\t r)});
	\draw (0.99,3.08) node[anchor=north west] {$\theta$};
	\draw [line width=1.2pt,dash pattern=on 1pt off 2pt on 5pt off 4pt] (2.12,0.06)-- (2.14,-2.29);
	\draw (2.84,-0.9) node[anchor=north west] {width$=r$};
	\begin{scriptsize}
	\draw [fill=qqqqff] (-1.56,2.75) circle (1.5pt);
	\end{scriptsize}
	\end{tikzpicture}
	\caption{A cut of a cone of angle $\theta$ and a cut of a cylinder of width $r$.}
	\label{cut}
\end{figure}
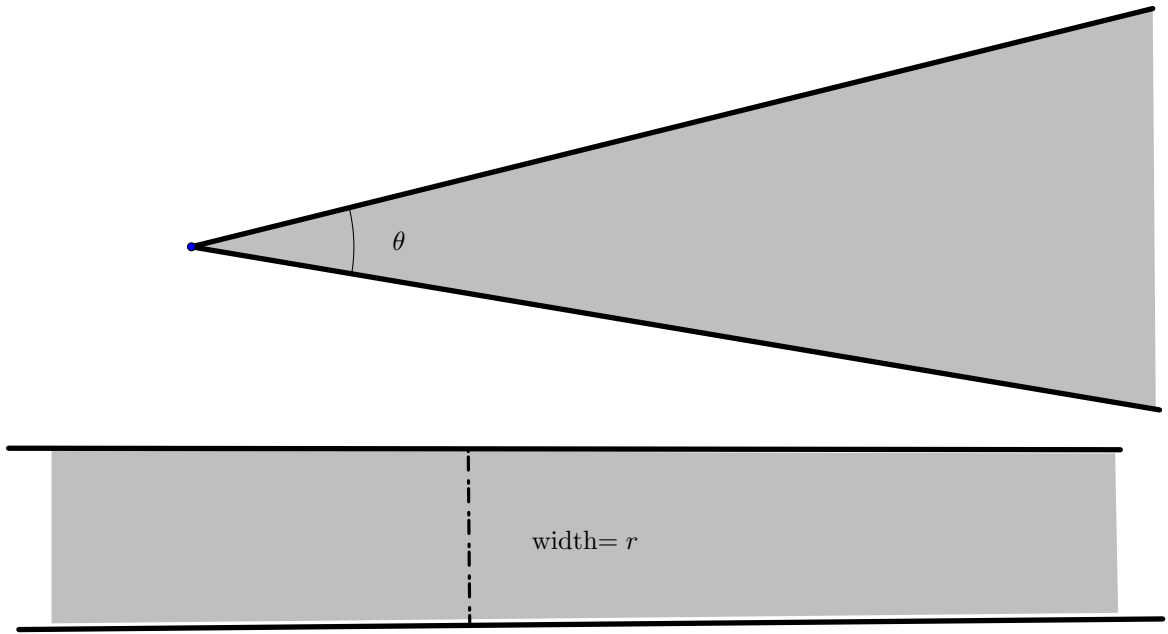

\begin{definition}
	A cut of a cone of angle $\theta$ is the once punctured  disk obtained by cutting a cone of angle $\theta$ along a geodesic directing from its origin, and will be denoted as $V_{\theta}$. 
\end{definition}
\noindent A cut of cone can be regarded as an FDL  disk with one punctured 
point and one labeled points at its boundary. As usual, angle and 
curvature at the labeled point are $\theta$ and $\pi-\theta$, respectively. For the punctured point,  angle and  curvature at 
the punctured point are $ -\theta$ and $\pi+\theta$, respectively.
Hence Gauss-Bonnet formula for the closed disk holds. See Figure \ref{cut}.

\begin{definition}
	A cut of a cylinder, $I(r)$, is the twice punctured disk together with a metric which is isometric to an infinite strip in the Euclidean plane. Width of the strip, $r$, is called width of the cut.
\end{definition}

\noindent A cut of cylinder can be regarded as an FDL  disk with two punctured points at its boundary. By definition, the angle and the curvature at each of the punctured  points are $0$ and $\pi$, respectively. Hence Gauss-Bonnet formula for the closed disks holds.
See the Figure \ref{cut}. 

\begin{definition}
	A cylinder of width $r$, $C_{0,r}$, is a metric space obtained by identifying edges of  a cut of a cylinder having  width $r$ through \textit{opposite} points.
\end{definition}

\noindent Observe that a cylinder can be considered as  FDL sphere with two punctured points. By convention, angles at these punctures are $0$. We can also call a cylinder as a cone of angle $0$. Also, again by convention, the curvature at each of the punctured points,  is $2\pi$. Observe that Gauss-Bonnet formula for the sphere 
holds.

\subsection{Cone metrics on  disk}

\begin{definition}
	A (flat) cone metric on a DL surface $S^L$ is a metric on $S_{\frak{p,p'}}$ so that each point $x$ in  $S_{\frak{p,p'}}$ has a neighborhood isometric to a neighborhood of the apex of a cone $C_{\theta}=C_{\theta_x}$ or a cut of a cone $V_{\theta}=V_{\theta_x}$, and
	\begin{itemize}
		\item
		$\frak{l}=\{y \in S_{\frak{p},B}: \ \theta_y \neq 2\pi \}$,
		\item
		$\frak{l'}= \{ y \in B-\frak{p'}: \ \theta_y \neq \pi   \}$.
	\end{itemize}
	Angle at $x$, $\theta(x)$ is defined to be $\theta_x$. If $x \in S_{\frak{p},B}$, then 
the	curvature at $x$, $\kappa(x)$, is defined as $2\pi - \theta(x) $. If $x \in B-\frak{p'}$, then the curvature is  $\kappa(x)=\pi -\theta(x)$. $x$ is called 
	singular if $\kappa(x)\neq 0$. Otherwise it is called non-singular.
	
\end{definition}
\noindent Observe that the two conditions above guarantee that set of singular points of $S^L$ is $\frak{l}\cup \frak{l'}$.

 Cones, cylinders, cuts of  cones, cuts of cylinders are examples of cone metrics. Observe that each FDL surface can be regarded as a cone metric on the underlying DL surface. Note that by an isometry of cone metrics on DL surfaces $S^L$ and $\bar{S}^L$, we mean an isometry of underlying metric spaces $S_{\frak{p,p'}}$ and 
$\bar{S}_{\frak{\bar{p},\bar{p'}}}$.
 
Now we state some elementary facts about cones, cylinders, cuts of cones and cuts of cylinders without proof.
\begin{proposition}
	\label{V}
	\label{proposition cylinder}
	\begin{enumerate}
		\item
		Let $d$ a complete  cone metric on a 1-punctured and 1-labeled DL  sphere  
		$S^L$. $S^L$ is isometric to $C_{\theta}$, for some $\theta>0$.

		\item 
		Let $d$ be a complete cone metric  on the  2-punctured DL sphere  $S^L$. $S^L$ is  isometric to $C_{0,r}$, for some $r>0$.
		
		\item

		Let $d$ be a complete  cone metric  on the  2-punctured DL  disk
		$S^L$, where the punctures are on the boundary. $S^L$ is isometric to $I(r)$, for some $r>0$.
		\item

		Let $d$ be a complete  cone metric  on the  a DL  disk $S^L$ with one punctured and one labeled boundary points. $S^L$ is isometric to a cut of a cylinder $V_{\theta}$, for some $\theta>0$. 
	\end{enumerate}
\end{proposition}

\begin{proposition}
	\begin{itemize}
		\item
		Two cones $C_{\theta}$ and $C_{\theta'}$ are isometric if and only if
		$\theta=\theta'$,
		\item
		Two cut of cones $V_{\theta}$ and $V_{\theta'}$ are isometric if and only if $\theta=\theta'$,
		\item
		Two cut of cylinders $I(r)$ and $I(r')$ are isometric if and only if
		$r=r'$,
		\item
		Two  cylinders $C_{0,r}$ and $C_{0,r'}$ are isometric if and only if $r=r'$.
	\end{itemize}
	
\end{proposition}


\subsubsection{ Cone metrics on  disk with one punctured and two labeled  boundary points}

\begin{figure}
	\begin{center}
		\includegraphics[scale=0.4]{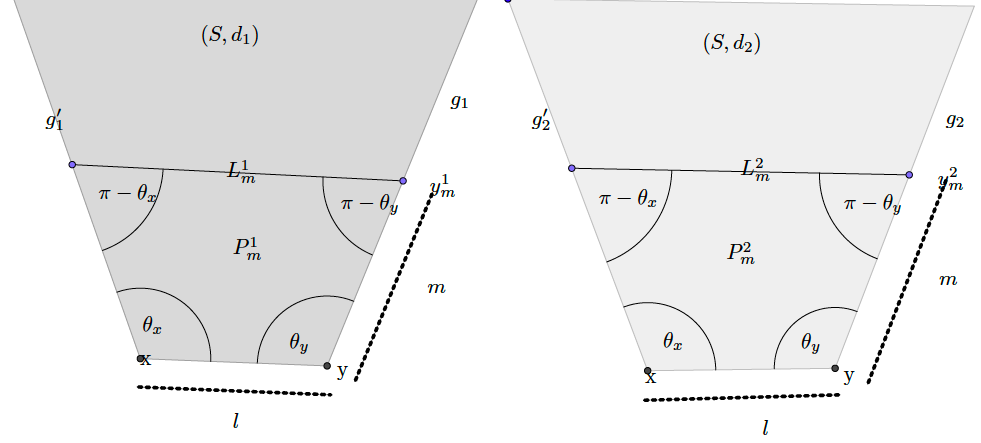}
	\end{center}
	\caption{}
	\label{disom}
\end{figure}

 Our next objective is classify cone metrics on a DL  disk $S^L$ with one punctured and 
2 labeled boundary points. See Figure \ref{disom}. Let $x,y$
be the labeled points of the boundary and $g$ be the part of the boundary which connects $x$ and $y$. Note that we want classify cone metrics up to isometries which fix $x$ and $y$.

\begin{proposition}
	\label{D}
	
	\begin{enumerate}
		\item
		For each  triple 
		of positive numbers $(\theta_x,\theta_y, l)$, 
		$\theta_x+\theta_y\geq\pi$, there exists a  complete cone metric on $S^L$ so that the angle at $x$ is 
		$\theta_x$, the angle at $y$ is $\theta_y$ and length of 
		$g$ is $l$.
		\item
		Each complete cone metric on $S^L$ is uniquely determined by its angles and length of $g$.
	\end{enumerate}
	\begin{proof}
		\begin{enumerate}
			\item
			\textup{
				There are two cases to be considered separately.}
			\begin{itemize}
				\item
				Assume one of $\theta_x$ and $\theta_y$ is less than or equal to $\frac{\pi}{2}$.
					Without loss of generality let this be $\theta_x$. Therefore 
					$\theta_y\geq \frac{\pi}{2}$ . If $\theta_x=\theta_y=\frac{\pi}{2}$, we know that there exists such a region in the Euclidean plane. If not, form a complete cone metric on a 2 labeled 1 punctured disk  with cone angles at $\theta_x$,  $\pi-\theta_x$ and length of the segment joining labeled points is $l$. Indeed, such a surface can be drawn in plane. Call the vertices of the surface with angles $\theta_x$ and $\pi-\theta_y$ as $x$ and $y$, respectively. Let's denote the half line on  the surface originating from $y$ as $L_y$.  
					Now, take a cut of a cone of angle $\theta_y-\pi+\theta_x$. Glue one of the boundaries of the cut of the cone with $L_y$.
				The resulting surface has the properties that we want.

				\item
				
					Assume $\theta_x,\theta_y>\frac{\pi}{2}$. By the first part, there exists a complete cone metric on a 2-labeled and 1-punctured disk with cone angles $\frac{\pi}{2}$, and $\theta_y$ and the length of the segment joining these labeled points is $l$. Call the vertices on this surface with angles $\frac{\pi}{2}$ and $\theta_y$ as $x$ and $y$, respectively. Let's denote half line on this surface originating from $x$ as $L_x$. Glue one of the boundary geodesics of the cut of a cone of angle $\theta_x - \frac{\pi}{2}$ with the $L_x$. Resulting metric has desired properties.
				
			\end{itemize}
			\item
			\textup{Take two complete cone metrics $d_1, d_2$ on $S^L$ with the same angle and length data, $(\theta_x,\theta_y,l)$. We will consider two cases separately:}
			
			\begin{itemize}
				
				\item 
				For each $i=1,2,$ let $g_i$ be the half line on the boundary of surface which is based at $y$, with respect to $d_i$. For each $i=1,2,$ let $g'_i$ be the half line on the boundary of surface which is based at $y$, with respect to $d_i$.
				 Let $y^{i}_m, m \in \mathbb{N},$ be the point on $g_i$ whose distance with  $y$ is $m$, with respect to $d_i$. Let 
				$L^{i}_m$ be the line segment joining $y^{i}_m$ with $g'_i$ so that the angle 
				between $L^{i}_m$ and the line segment $[y,y^{i}_{m}]$ is $\pi-\theta_y$.
				If we cut $(S,d_i)$ through $L^i_m$,  we will get convex polygons, $P^i_m$, for each $i=1,2$ and for 
				each $m\in \mathbb{N}$, which are evidently isometric.  
				For each $i=1,2 $, $\cup_{m\in \mathbb{N}}P^i_m=S_{\frak{p,p'}}$, therefore $(S, d_1)$ and $(S,d_2)$ are isometric. See Figure \ref{disom}. 
				\item 
				If one of  $\theta_x$ and $\theta_y$ is greater than or equal to $\pi$,
				one can cut both of the cone metrics through half-lines originating 
				from $x$ and $y$ to reduce the problem to the previous case. We omit the  details.
				

		\end{itemize}
			
		\end{enumerate}

	\end{proof}

\end{proposition}

  $S^L$ together with such a cone metric will be denoted as  $D(\theta_1,\theta_2, l)$.

\begin{remark}
	There is no complete cone metric on $S^L$ having angle data $(\theta_1,\theta_2)$ so that $\theta_1+\theta_2< \pi$.  
\end{remark}

\begin{remark}
	Assume that $\theta_1+\theta_2\geq \pi$. For each positive real number $r$, $D_{\bar{\kappa}}(l)$ can be triangulated so that
	\begin{enumerate}
		\item
		The length of edges of triangles lying in half-lines of the boundary is $r$.
		\item
		The triangulation satisfies properties in Definition \ref{triangulation}. 
		\item
		The metric obtained by  triangulation is the exactly that of $D_{\bar{\kappa}}(l)$.
	\end{enumerate}
	
	One can manage  to do this by decomposing $D_{\bar{\kappa}}(l)$ as in the proof of Proposition \ref{D}.
	\label{DT}
\end{remark}

\subsection{Cone metrics on the closed disk with one punctured interior or one punctured boundary point}

Let $D_{1,n}^L$ be a DL  disk with one punctured point at its interior
and $n$ labeled points on its boundary so that $\mathfrak{p'}$ and $\mathfrak{l}$  are empty. Similarly, let $\bar{D}^L_{1,n}$ be a DL disk with one punctured  and $n$ labeled boundary points so that $\mathfrak{p}, \mathfrak{l}$ are empty.
The aim of this section is to give  a complete classification of cone metrics of non-positive curvature on $D_{1,n}^L$. It turns out that  the length and the curvature data on  the boundary of such a disk explicitly describe the cone metric. We also give a similar result for the case of $\bar{D}^L_{1,n}$.

\begin{lemma}
	Consider a complete cone metric on $D_{1,n}^L$ and a boundary point $x$. Assume that curvature at each boundary 
	point is  non-positive. Let $g$ be a geodesic starting at $x$ and pointing the interior of $D_{1,n}^L$. $g$ does not hit the boundary and it is not self intersecting.
	
	\begin{proof}
		If curvature at each boundary point is $0$, then $D_{1,n}^L$ is isometric to half of a cylinder and the statement is true. Assume that this is not the case. Observe  that a geodesic with above properties can not intersect itself without winding once around 
		the puncture. Otherwise, we get a disk with only one singular point, and this singular point is on the boundary. Clearly, such a disk can not exist. Assume that
		it intersects boundary or itself. 
		
		There are two cases to be considered. First, consider the case 
		in which the geodesic intersects  the boundary \textit{before} it intersects itself. In that case, some part of the geodesic and the boundary form a polygon which has at most two vertices, having angle less than $\pi$. Vertices at intersection of the geodesic with the boundary,  and all the other
		vertices have angle bigger than or equal to $\pi$. By Gauss-Bonnet theorem for compact surfaces with boundary
		\cite{Tro-compact}, such a polygon does not exist.
		
		Second, assume that  the geodesic \textit{first} intersects itself. We can cut $S^L$ through the loop formed by the geodesic and obtain a cone metric on a closed annulus. Total curvature for the boundary component of the annulus,
		which results from the boundary of $D_{1,n}^L$, is negative. Total curvature
		for the other boundary component is non-positive. Indeed, it contains at most one singular point which has non-positive curvature. This contradicts with Gauss-Bonnet theorem \cite{Tro-compact} since such an annulus should have \textit{zero} total curvature.
	\end{proof}
\end{lemma}

We point out that geodesics on disks above tend to the punctured point, or the \textit{point at infinity of the disk}.
Let $b_1, b_2, \dots b_n$ be the labeled points  given in a cyclic order on the boundary. Let
$\kappa_i$ and $l_i,\ i=1,\dots,n,$ be real numbers so that $\kappa_i < 0$ and $l_i>0$ for each 
$i$.

\begin{figure}
	\begin{center}
		\includegraphics[scale=0.4]{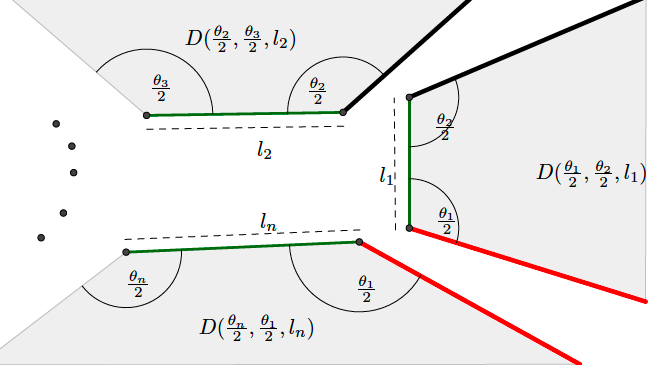}
	\end{center}
	\caption{Complete cone metrics on once punctured disk can be obtained from cone metrics on the disk with one puncture on its boundary.}
	\label{Ds}
\end{figure}

\begin{lemma}
	
	There exists a complete cone metric on $D_{1,n}^L$ so that
	for each $i=1,2 \dots, n$ the curvature at $b_i$ is $\kappa_i$
	and the length of $[b_i,b_{i+1}]$ is $l_i$.
	\begin{proof}
		\textup{Let $\theta_i=\pi - \kappa_i$. Consider $D(\frac{\theta_1}{2},\frac{\theta_2}{2},l_1), \dots, D(\frac{\theta_n}{2},\frac{\theta_1}{2},l_n) $. For $i<n-1$,
			glue $D(\frac{\theta_i}{2},\frac{\theta_{i+1}}{2},l_{i})$, along the geodesic 
			originating from the vertex having angle ${\frac{\theta_{i+1}}{2}}$ , with
			$D(\frac{\theta_{i+1}}{2},\frac{\theta_{i+2}}{2},l_{i+1})$, along the geodesic originating from the vertex having angle $\frac{\theta_{i+1}}{2}$. Do the same for $D(\frac{\theta_n}{2},\frac{\theta_1}{2},l_n) $ and
			$D(\frac{\theta_1}{2},\frac{\theta_2}{2},l_1)$. One will get a metric of desired type. See Figure \ref{Ds} }
	\end{proof}
	
	\label{ext}
\end{lemma}

\begin{remark}
	Assume that $D_{1,n}^L$ has a complete metric with curvature at each boundary point less than 
	or equal to $0$.  Let $x$ and $y$ be two distinct boundary points so that there is a straight boundary segment joining them. Take 
	two half-lines originating from $x$ and $y$ which are perpendicular to the segment considered. Gauss -Bonnet theorem implies that these two half-lines does not intersect.
\end{remark}

\begin{lemma}
	Consider $D_{1,1}^L$ . Let $\kappa_1$ and $l_1$ be as above. Then there is a unique cone metric on $D_{1,1}^L$ having curvature $\kappa_1$ at
	$b_1$ and the length of the boundary is $l_1$.
	\begin{proof}
		We proved existence of the metric. See Lemma \ref{ext}. Take such a metric on $D^L_{1,1}$. Let $\theta=\pi-\kappa_1$. If we cut $D_{1,1}^L$ through the geodesic making an angle $\frac{\theta}{2}$ with the boundary, then the resulting surface. Therefore it is isometric to $D(\frac{\theta}{2},\frac{\theta}{2},l_1)$. If we glue back, we get the metric we started. Thus, any metric with these properties obtained by gluing the half lines of the boundary of 
		$D(\frac{\theta}{2},\frac{\theta}{2}, l)$.
		Hence, there exist a unique metric having properties stated in the present lemma.
	\end{proof}
	
	\label{1S}
\end{lemma}

\begin{theorem}
	\label{class1}
	Let $\kappa_i<0$ and $l_i>0$, $i=1,\dots,n$, be real numbers. There is a unique complete cone metric on $D_{1,n}^L$, up to isometries respecting labeling, so that curvature at $b_i$ is $\kappa_i$
	and length of the segment $[b_i,b_{i+1}]$ is $l_i$ for each $i=1,\dots,n$. 
	\begin{proof}\textup{
			We proved existence of such a metric. See Lemma \ref{ext}.} 
		\par{\textit{Uniquness}}\\
		\textup{We use induction on  number of labeled points to prove the statement. Lemma \ref{1S} asserts that the statement is true  if number of  labeled points is one. Assume that the statement is true for the case that there are 
			$n$ or less labeled points. Let $d_1$ and $d_2$ be  metrics on $D_{1,n+1}^L$ having same curvature data.  Consider the the segment joining $b_n $ and $b_{n+1}$, call it $g$. By assumption, 
			$g$ has the same length with respect to two metrics. For each $i=1,2$, let 
			$g_i$ and $h_i$ be the half-lines originating from $b_n$ and $b_{n+1}$, with respect to $d_i$, so that $g_i$ and $h_i$ are perpendicular to $g$. Cut $(D_{1,n}^L,{d_i})$ through $g_i$ and $h_i$. For each $i$
			we get two DL surfaces $S_i$ and $D(\frac{\pi}{2},\frac{\pi}{2},l)$ where $l$
			is the length of the segment $g$. Glue $S_i$ through the half-lines on the boundary to get complete cone metrics on the disk with $n$ labeled points and one puncture. Call these surfaces, together with induced metrics, $(S'_i,d'_i)$. By induction hypothesis $(S'_1,d'_1)$ and $(S'_2,{d'_2})$ are isometric. Thus, 
			 metrics on $D_{1,n}^L$ obtained from $d_1'$ and $d_2'$ by reversing the cutting and gluing operation above   are same. Therefore, these induced metrics  should coincide with $d_1$ and $d_2$.
			Hence $(D_{1,n}^L,d_1)$ and $(D_{1,n}^L,d_2)$ are isometric.}
	\end{proof}
	\label{theo}
\end{theorem}
\begin{remark}
	\label{sumatboundary}
	If a cone metric on $D^L_{1,n}$ is complete, then $\sum_{x \in \frak{l'}}\kappa(x)\leq 0$.
\end{remark}
A DL surface together with the  metric having curvature data
$\bar{\kappa}=(\kappa_1,\dots, \kappa_n)$ and length data $\bar{l}=(l_1,\dots, l_n)$ will be denoted as $D_{\bar{\kappa}}(\bar{l})$, where $\kappa_i < 0, l_i>0$.

\begin{corollary}
	\label{tribar}
	
	$D_{\bar{\kappa}}(\bar{l})$  can be triangulated so that the triangulation 
	has properties in Definition \ref{triangulation} and the induced metric coincides 
	with the metric of $D_{\bar{\kappa}}(\bar{l})$.
	\begin{proof}
		\textup{By  Theorem \ref{theo}, $D_{1,n}^L$ can be decomposed into finite numbers of disks of the form $D(\theta_1,\theta_2,l)$. Hence the result follows from Remark \ref{DT}.}
	\end{proof}
\end{corollary}

\begin{corollary}
	Assume that $\kappa_i$ and $l_i$ satisfies the above conditions , and also $\kappa_i=\kappa_j=\kappa>-\pi$ and $\l_i=\l_j=l$ for all $i,j=1,\dots,n$.  $D_{\bar{\kappa}}(\bar{l})$ can be embedded in a cone.
	\begin{proof}\textup{
			Consider the cone with  angle 
			$-n\kappa$. Obviously, there is a compact polygonal part of the cone, homeomorphic  to a disk, having the apex as an interior point and $n$ boundary edges of length $l$, $n$ boundary points of angle $\pi+\kappa$. Closure of the complement of this disk has the same length and curvature data with that of $D_{\bar{\kappa}}(\bar{l})$. The result follows from the uniqueness part of the above theorem.}
	\end{proof}
\end{corollary}

\begin{figure}
	\begin{center}
		\includegraphics[scale=0.4]{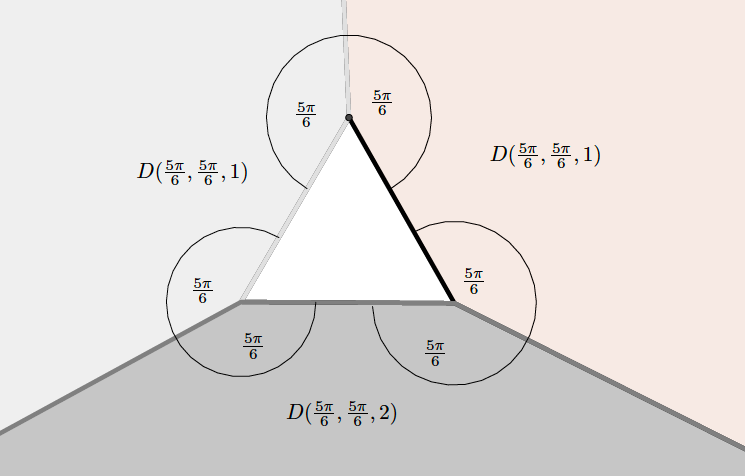}
	\end{center}
	\caption{A non-planar cone metric on once punctured disk. This disk also has an irregular puncture.}
	\label{nonplane}
\end{figure}

\begin{example}
	Consider the cone metric on $D_{1,3}^L$ obtained by  gluing 
	two copies of $D(\frac{5\pi}{6},\frac{5\pi}{6},1)$ and one copy of $D(\frac{5\pi}{6},\frac{5\pi}{6},2)$. $D^L_{1,3}$, together with this metric, can not be embedded into a cone. Otherwise, by Gauss-Bonnet Formula   this cone would have  angle at its apex equal to $2\pi$. Hence, it would be  the Euclidean plane. This embedding produces a triangle in the plane  having edge lengths $1,1,2$ and angles $\frac{\pi}{3},\frac{\pi}{3},\frac{\pi}{3}$ on the plane, which does not exist. Also observe that 
	one can not embed $D_{1,3}^L$ into a cone even after removing any compact set. This means that $D^L_{1,3}$ has an irregular puncture. See Figure \ref{nonplane}.
\end{example}

Now we state analogous results for $\bar{D}_{1,n}^L$. We omit the proofs since they are entirely analogous to the proofs of the facts we obtained for complete cone metrics 
on $D_{1,n}^L$. Assume that $b_1, \dots b_n$ are labeled boundary points 
so that $b_1, \dots, b_n$ and the punctured point are in a cyclic order in boundary. Note that this labeling implies that $b_1$ and $b_n$  share same edges with the puncture.

\begin{theorem}
	\label{dbp}
	Assume that we are given numbers $\kappa_1,\dots, \kappa_n$,  $n\geq 1$, so that 
	\begin{itemize}
		\item
		$\kappa_1,\kappa_n<\pi$,
		\item
		$\kappa_2,\dots,\kappa_{n-1} < 0$
		\item
		$\sum^n_{i=1}\kappa_i\leq\pi$
	\end{itemize}
	and $\l_1,\dots,\l_{n-1}$ so that $l_i>0$ for each $i$. There exists a unique complete cone metric on $\bar{D}_{1,n}^L$ so that curvature at $b_i$ is $\kappa_i$ and length of the segment $[b_i,b_{i+1}]$ is $l_i$. Also, two cone metrics having same length and curvature data are isometric.

\end{theorem}

\begin{remark}
	\label{sumatboundary}
	If a cone metric on $\bar{D}^L_{1,n}$ is complete, then $\sum_{x \in \frak{l'}}\kappa(x)\leq \pi$.
\end{remark}
	
We will denote $\bar{D}_{1,n}^L$ together with such a metric as $\bar{D}_{\bar{\kappa}},(\bar{l})$,
where $\bar{\kappa}=(\kappa_1,\dots,\kappa_n)$, $\bar{l}=(l_1,\dots,l_{n-1})$.
Note $\bar{D}^L_{1,1}$ is nothing else than a cut of cone $V_{\theta}$. See Proposition \ref{V}. 
\begin{corollary}
	\label{trinobar}
	Let  $\bar{\kappa}$, $\bar{l}$ be  the curvature and the length data satisfying properties in Theorem \ref{dbp}.  $\bar{D}^L_{1,n}$ can be triangulated so  that the triangulation satisfies the properties of Definition \ref{triangulation} and the induced cone metric coincides with the metric of $D_{\bar{\kappa}}(\bar{l})$.
	
\end{corollary}

\subsection{Modification }
\label{modification1}
We conclude this section by some results about cone metrics on DL  disks without labeled interior points.





\paragraph{Modification:}
\label{modification-definition}
By a modification of a cone metric on a DL closed disk  without  labeled interior points, we mean the resulting surface (with the induced metric)
obtained after recursively cutting finitely many Euclidean triangles which are incident to the boundary at least at one edge. Note that we require the triangles shold be incident with the boundary at most at one edge.

\begin{figure}
	\begin{center}
		\includegraphics[scale=0.50]{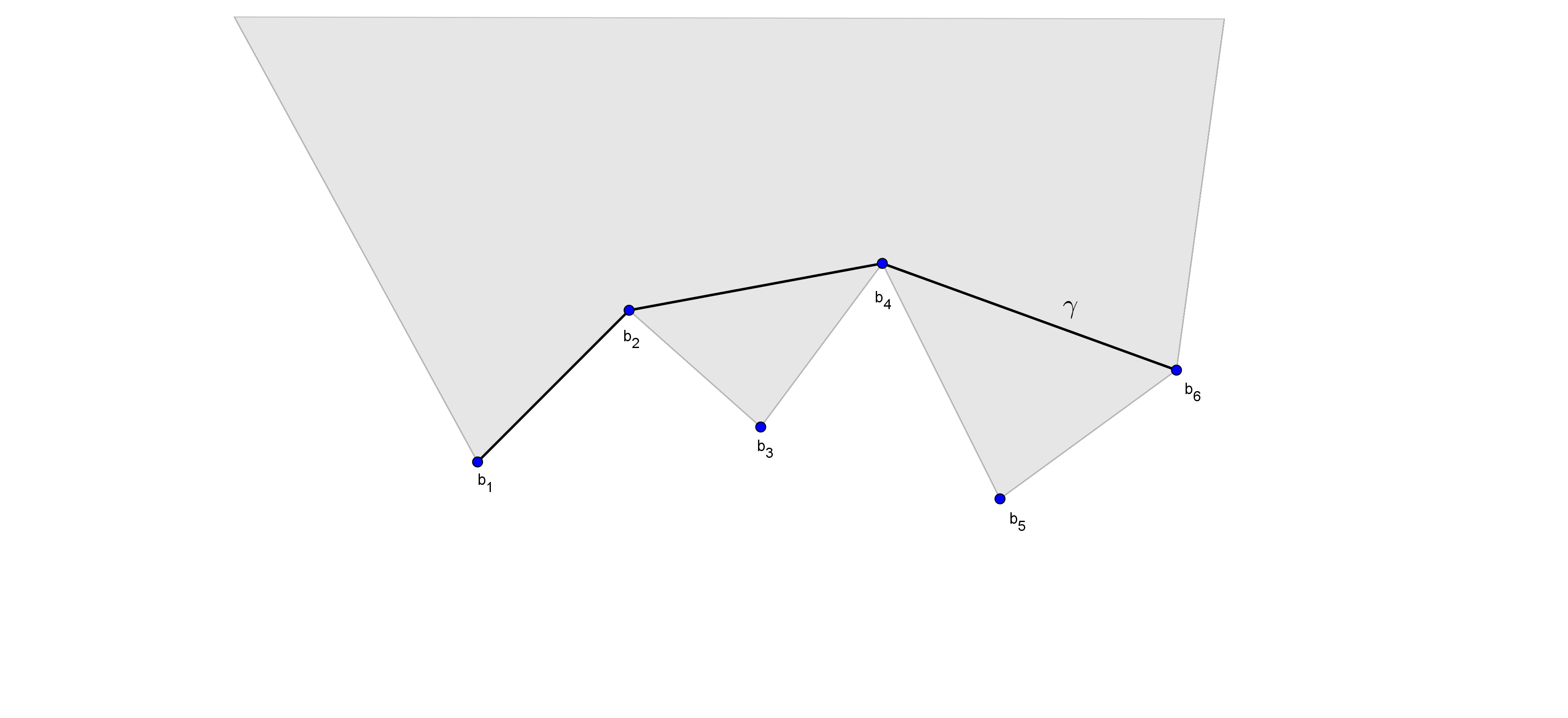}
	\end{center}
	\caption{ $\bar{D}^L_{1,m}$ and $\gamma$. $\gamma$ is length minimizing curve joining $b_1$ with $b_m$.}
	\label{gamma}
\end{figure}

\begin{figure}
	\begin{center}
		\includegraphics[scale=0.50]{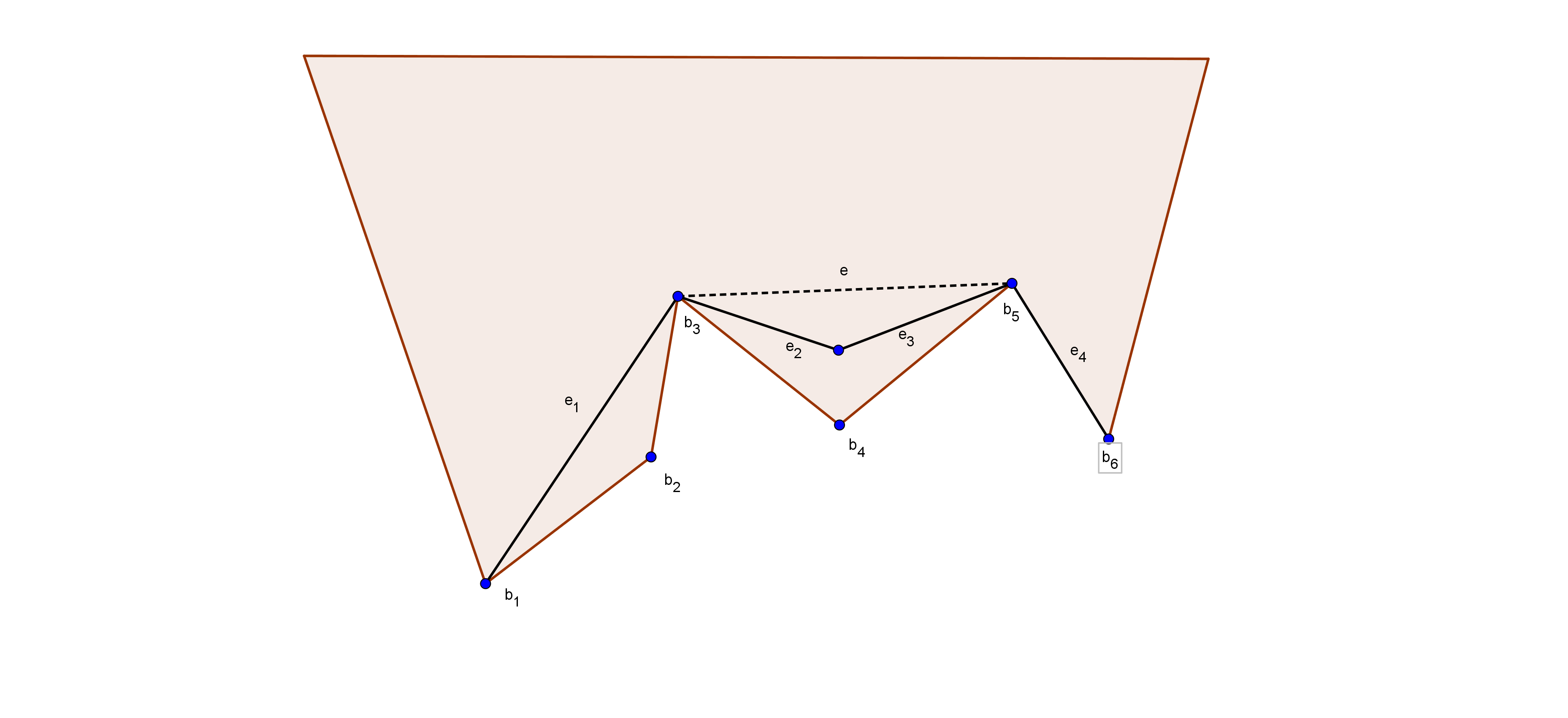}
	\end{center}
	\caption{The curve with edges $e_1,\ e_2,\ e_3,\ e_4$ can not be length minimizing since $\lvert e_2\rvert + \lvert e_3\rvert> \lvert e\rvert$. Note that we denote length of an edge $E$ by $\lvert E\rvert$.} 
	\label{notgamma}
\end{figure}

\begin{proposition}
	\label{mod1}
	Every complete cone metric on $\bar{D}_{1,n}^L$ can be modified as follows:
	\begin{enumerate}
		\item
		If $x \in B$, then $\theta(x)<2\pi$.
	\item
	If $x \in \frak{l'}$ and not sharing an edge with the puncture, then $\kappa(x)<0$. 
	\end{enumerate}
	\begin{proof}
		First of all, consider $\bar{D}^L_{1,1}$. A complete cone metric on it is nothing
		else than a cut of cone, and the statement is true if its angle $\theta<2\pi$. The statement is also is true
		for $\bar{D}^L_{1,0}$ which is half plane.
		
		Consider a complete flat metric on $\bar{D}^L_{1,n}$, $n\geq 2$, or on $\bar{D}^L_{1,1}$ where angle at the singular vertex is greater than or equal to  $2\pi$.
		\begin{itemize}\item
		  If there is a boundary point having an angle greater than or equal to $2\pi$, then we can remove a polygon about it so that resulting singular points have angle less than $2\pi$.  Therefore, by removing a polygon about each singular point $x$ such that $\theta(x)\geq 2\pi$, we get a complete DL disk for which boundary 
		  points have angle less than $2\pi$.
		  \item
 If after above operation we get $\bar{D}^L_{1,0},\ \bar{D}^L_{1,1}$ or $\bar{D}^{L}_{1,2}$, then we are done.
 Thus assume that we get a complete cone metric on $\bar{D}^{L}_{1,m}, \ m\geq 3$. Let's label its singular points as $b_1, \dots b_m$. See
 Figure \ref{gamma}. Take a loop joining $b_1$ with
 $b_n$. There is a length minimizing curve in its homotopy class. See \cite{Gromov}. Call this $\gamma$. If we cut $D^L_{1,m}$, then we get a surface of the type we want. If this is not 
 the case, then $\gamma$ has two edges making
 angle less than $\pi$. This implies that $\gamma$ is not length minimizing. See Figure  \ref{notgamma}. 
		 \end{itemize}
	\end{proof}
\end{proposition}

\begin{figure}
	\begin{center}
		\includegraphics[scale=0.5]{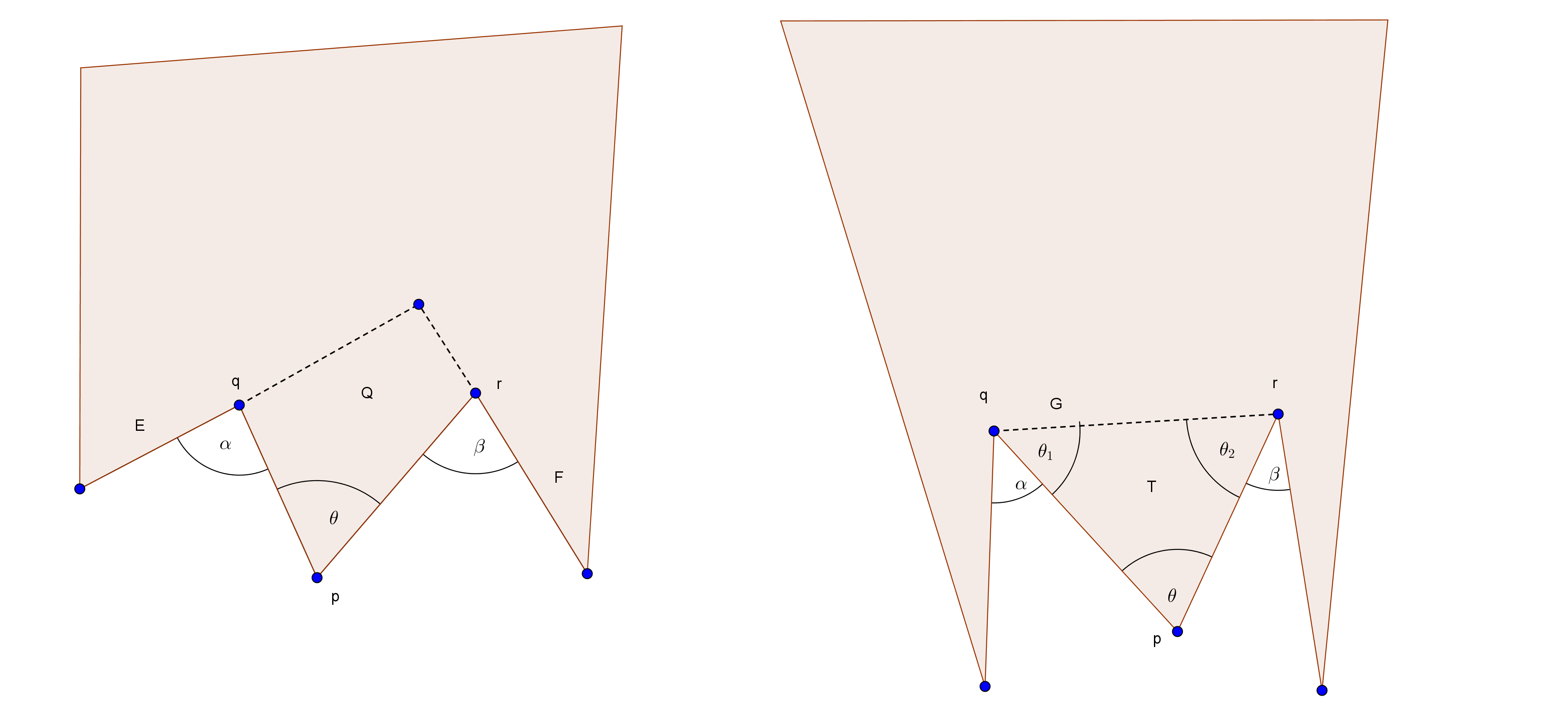}
	\end{center}
	\caption{
			 If $\theta <\alpha+\beta$, then we  remove the quadrilateral $Q$ to obtain desired modification. See left of the figure.
			If $\theta\geq \alpha+\beta$, then we remove the triangle $T$ to obtain desired modification. See right of the figure.}
	\label{bilemedim}
\end{figure}

There is a similar result for the cone metrics on the closed disk having one punctured interior point. The proof is also similar. We state it and outline its proof.
\begin{lemma}
	\label{mod2}
	Each complete cone metric on a $D_{1,n}^L, n\geq2$   can be modified so that resulting surface has the following properties:
	
	\begin{enumerate}
		\item 
		Each point on the boundary has angle less
		than $2\pi$.
		\item
		There is at most one singular point of positive curvature.
	\end{enumerate}
	
	\begin{proof}
		First modify the cone metric so that there are 
		no singular points having angle greater than and equal to $2\pi$. One can do this 
		as in the proof of Proposition \ref{gamma}.
		Then take a boundary point $x$ and consider a loop based at $x$ and winding once around the puncture. Take a length minimizing path in the homotopy class of the loop and cut the surface through this path. Resulting cone metric has at most one singular point of
		positive curvature, $x$.
	\end{proof}
	
\end{lemma}
\begin{proposition}
	\label{mod3}
	Let $D^L_{1,n}$ be an FDL surface so that 
	\begin{enumerate}
		\item 
		It has one singular point of positive curvature,
		\item
		Each boundary point has angle less than $2\pi$.
	\end{enumerate}
$D^L_{1,n}$  has a modification so that for each 
boundary point $x$, $\pi\leq \theta(x) < 2\pi$.
\begin{proof}
	Since total curvature at the boundary of  $D^L_{1,n}$ is non-positive, $n\neq 0,1$. Assume that $n=2$. Let $p$ be the singular point with positive curvature. Take a length
	minimizing loop which is based at $p$ and winds once around boundary. If we cut the disk through this loop, we get a disk at most one 
	singular point, $p$. Curvature at $p$ is not positive, since modification does no change 
	total curvature. Also, it is clear that the angle at $p$ is less than $2\pi$.
	
	We do induction on number of singular points.  
	 Consider a flat metric on $D^L_{1,n}$, $n\geq3$. We denote the singular point with positive curvature by $p$, and singular points which share an $edge$ with $p$ by $q$ and $r$.  See Figure \ref{bilemedim}. In that figure $\theta$ is the angle at $p$, $\alpha=\pi+\kappa(q)$, $\beta=\pi+\kappa(r)$.
	 There are two cases to be considered:
	 \begin{enumerate}[a)]
	 	\item
	 	$\theta<\alpha+\beta$. In this case, we can extend the edges $E$ and $F$ to form the quadrangle $Q$. See left of the Figure \ref{bilemedim}. If we remove $Q$, we get a surface with at most one singular point of positive curvature, and it is clear that the number of singular points of these surface is less than $n$. 
	 	\item
	 	
	 $\theta\geq \alpha+\beta$. Draw a line segment joining $q$ and  $r$ to form a triangle. Call the segment $G$ and the triangle $T$. See right of the Figure \ref{bilemedim}. Let $\theta_1$ and $\theta_2$
	 be the angles at $q$ and $r$, respectively.
	 Since $\pi-\theta_1-\theta_2=\theta$, we have
	 $$\pi-\theta_1-\theta_2\geq \alpha+\beta$$
	 $$\pi\geq \alpha+\theta_1+\beta+\theta_2.$$
	 
	 	 Therefore one of $\alpha+\theta_1$ and $\beta+\theta_2$ less than $\pi$. This means that when we remove the triangle we reduce number of singular points and the resulting surface has at most one singular points of positive curvature.	
	 \end{enumerate}
\end{proof}
\end{proposition}

\begin{corollary}
	\label{mod-corollary}
	Each complete cone metric on $D^L_{1,n}$, $n\geq 0$, can be modified so that resulting disk
	does not have points with positive curvature on its boundary.
	\begin{proof}
		The statement immediately follows from 
		Lemma \ref{mod2} and Proposition \ref{mod3}
	\end{proof}
\end{corollary}

\section{DL surfaces with complete cone metrics can be triangulated}
\label{triangulation-subsection}
In this section we prove that DL surfaces together with complete cone metrics can be triangulated  so that resulting metric coincides with the given one. This theorem is well-known for compact surfaces \cite{Tro-compact}. Our strategy is to cut such a surface around its punctures and reduce the problem to the cases for compact surfaces and the surfaces $\bar{D}(\bar{\kappa},\bar{l}),D(\bar{\kappa},\bar{l})$.

\begin{theorem}
	\label{theorem-triangle}
	Every DL surface $S^L$ together with a complete metric $d$ can be triangulated as in Definition \ref{triangulation} so that resulting metric coincides with $d$.
	
	\begin{proof}
		For each point $p \in \mathfrak{p}$, take a non-self intersecting polygonal loop around $p$ so that the punctured disk bounded by $p$ and the loop has no labeled point on its interior and no punctured point on it except $p$. For each point $p' \in \mathfrak{p'}$, take a polygonal path joining to half-lines \textit{incident} $p'$ so that the punctured disk
		bounded by $p'$ and the path has no labeled points on its interior and no punctured point 
		on it except $p'$. Also observe that we can choose these loops and paths so that resulting disks are pairwise disjoint. Note that we may assume these disks satisfy properties in Proposition \ref{mod1} and Corollary
		 \ref{mod-corollary}.
		Now we know that these disks can be triangulated nicely. See Definition\ref{triangulation}  and Corollaries \ref{tribar}, \ \ref{trinobar}.  If we remove interiors of  these disks what we get is a compact surface together with a cone metric. It is well known that such a surface can be triangulated with only finitely many triangles. 
		Therefore we can use triangulations on these pieces to obtain a triangulation on $S^L$.  This triangulation has the properties in Definition \ref{triangulation} and  metric $d$ coincides with the metric induced by the triangulation at each triangle. Hence these metrics coincide globally.

	\end{proof}
	
\end{theorem}


\section{Gauss-Bonnet formula}
\label{Gauss-Bonnet}
Gauss-Bonnet formula for compact flat surface is well-known. There is a variant of the formula for the non-compact case. But it pre-assumes that 
each punctured interior point has a neighborhood isometric to a neighborhood of point at infinity of a cone. See \cite{Tro-open}, \cite{Tro-compact}.

We start with defining curvature at the punctures of a DL surface with a complete cone metric. Then we will state and prove Gauss-Bonnet Formula. 

\begin{remark}
	\label{modremark}
A modification of a complete cone metric on  $D^L_{1,n}$ does not change total  curvature of the boundary of $D^{L}_{1,n}$. 
Similarly, a modification of a complete cone metric on  $\bar{D}^L_{1,n}$ does not change total  curvature of the boundary of $\bar{D}^{L}_{1,n}$.
\end{remark}

\begin{definition}
	\begin{itemize}
		\item
		If $\bar{D}^{L}_{1,n}$ has a complete flat metric, then the  curvature at its puncture $p'$
		is defined as 
		$$\kappa(p')=2\pi- \sum_{x \in \mathfrak{l'}}\kappa(x).$$ 
		The angle at $p'$ is $\theta(p')=\pi-\kappa(p')$.
		\item
			If $D^{L}_{1,n}$ has a complete flat metric, then the curvature at its puncture $p$
		is defined as
		$$\kappa(p)=2\pi- \sum_{x \in \mathfrak{l'}}\kappa(x).$$ 
		The angle at $p$ is $\theta(p)=2\pi-\kappa(p)$.
		\item
		Let $S^L$ be a DL surface (together with a complete cone metric) and $p\in \frak{p}$. The curvature at $p$, $\kappa(p)$, is the curvature of $p$ as a
		punctured point of a disk in $S^L$ containing $p$ and having no singular points on its interior. The angle at $p$ is $\theta(p)=2\pi-\kappa(p)$.
		
		\item
		Let $S'^L$ be a DL surface (together with a complete cone metric) and $p'\in \frak{p'}$. The curvature at $p'$, $\kappa(p')$, is the curvature of $p'$ as a
		punctured point of a disk in $S'^L$ containing $p'$ and having no singular points on its interior. The angle at $p'$ is $\theta(p')=\pi-\kappa(p')$.

	\end{itemize}
\end{definition}

\begin{remark}
	By Remark \ref{modremark}, last two items of the above definition make sense. Any two such disks containing $p$ can be modified to a common disk, hence have same total curvature at their boundaries.
\end{remark}

\begin{theorem}[Gauss-Bonnet formula]
	Let $S^L$ be a DL surface together with a complete cone metric. The following formula holds: 
	
	\begin{align}
	\sum_{x\in S}\kappa(x)= 2\pi\chi(S).
	\end{align}
	
	\begin{proof}
		Assume that $S^L$ has $n$ punctured points on its interior and $m$ punctured points on its boundary. 
		As in the proof of Theorem \ref{theorem-triangle}, choose  \textit{disks} around the punctures. Let $S'$ be the compact surface, with induced metric, obtained by removing these disks. Observe that
		
		\begin{itemize}
			\item
			$\chi(S')=\chi(S)-n$, and
			\item
			$\sum_{y \in S'}\kappa(y)=2\pi\chi(S')$
		\end{itemize}
		by Gauss-Bonnet Formula for compact surfaces. Now, observe that removing one appropriate disk around a punctured interior point decreases total curvature $2\pi$. Therefore, if we remove $n$ such disks, total 
		curvature decreases $2n\pi$. Also observe that removing an appropriate disk around a punctured boundary point does not change total curvature. Therefore we have 
		
		$$\sum_{x \in S}\kappa(x)= 2n\pi+\sum_{y \in S'}\kappa(y)=
		2n\pi+2\pi\chi(S')= 2\pi \chi(S).$$

	\end{proof}
	
\end{theorem}

\section{Existence of geodesic representatives in free homotopy classes of loops}
\label{free-homotopy}

It is well-known that any loop in any compact flat surface  has a length minimizing closed geodesic representative in its \textit{free} homotopy class \cite{Gromov}.	Recall that by a length minimizing closed geodesic we mean a closed geodesic which 
has length less than or equal to length of each curve in its free homotopy class. We prove this property is valid  for any FDL surfaces. The idea of our proof is to cut the surface through
the punctures and reduce the problem to the case of compact surfaces. We start with some observations.

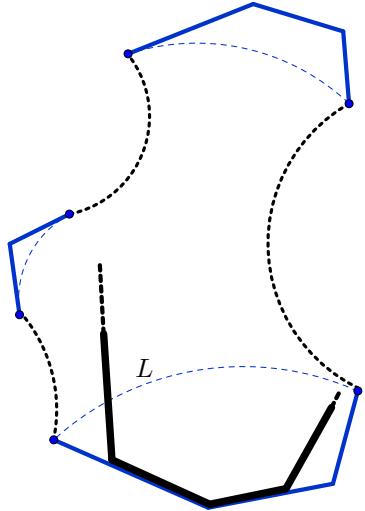
\begin{figure}
	\definecolor{qqttcc}{rgb}{0,0.2,0.8}
	\definecolor{qqqqff}{rgb}{0,0,1}
	\begin{center}
		\begin{tikzpicture}[line cap=round,line join=round,>=triangle 45,x=0.65cm,y=0.65cm]
		\draw [shift={(0.18,3.52)},line width=1.2pt,dotted]  plot[domain=-1.37:0.67,variable=\t]({1*2.04*cos(\t r)+0*2.04*sin(\t r)},{0*2.04*cos(\t r)+1*2.04*sin(\t r)});
		\draw [shift={(7.94,0.92)},line width=1.2pt,dotted]  plot[domain=2.09:4.26,variable=\t]({1*3.3*cos(\t r)+0*3.3*sin(\t r)},{0*3.3*cos(\t r)+1*3.3*sin(\t r)});
		\draw [shift={(-2.72,-2.48)},line width=1.2pt,dotted]  plot[domain=-0.21:0.71,variable=\t]({1*3.04*cos(\t r)+0*3.04*sin(\t r)},{0*3.04*cos(\t r)+1*3.04*sin(\t r)});
		\draw [line width=1.6pt,color=qqttcc] (0.26,-3.1)-- (3.42,-4.49);
		\draw [line width=1.6pt,color=qqttcc] (3.42,-4.49)-- (5.97,-4);
		\draw [line width=1.6pt,color=qqttcc] (5.97,-4)-- (6.48,-2.1);
		\draw [line width=1.6pt,color=qqttcc] (1.78,4.8)-- (4.34,5.82);
		\draw [line width=1.6pt,color=qqttcc] (4.34,5.82)-- (6.18,5.26);
		\draw [line width=1.6pt,color=qqttcc] (6.18,5.26)-- (6.3,3.78);
		\draw [line width=1.6pt,color=qqttcc] (0.58,1.52)-- (-0.64,0.91);
		\draw [line width=1.6pt,color=qqttcc] (-0.64,0.91)-- (-0.44,-0.54);
		\draw [shift={(4.13,-7.35)},dash pattern=on 2pt off 2pt,color=qqttcc]  plot[domain=1.15:2.31,variable=\t]({1*5.75*cos(\t r)+0*5.75*sin(\t r)},{0*5.75*cos(\t r)+1*5.75*sin(\t r)});
		\draw [shift={(3.18,0.46)},dash pattern=on 2pt off 2pt,color=qqttcc]  plot[domain=0.82:1.88,variable=\t]({1*4.56*cos(\t r)+0*4.56*sin(\t r)},{0*4.56*cos(\t r)+1*4.56*sin(\t r)});
		\draw [shift={(1.52,-0.23)},dash pattern=on 2pt off 2pt,color=qqttcc]  plot[domain=2.06:3.3,variable=\t]({1*1.98*cos(\t r)+0*1.98*sin(\t r)},{0*1.98*cos(\t r)+1*1.98*sin(\t r)});
		\draw [line width=2.8pt] (1.28,-0.93)-- (1.45,-3.51);
		\draw [line width=2.8pt] (1.45,-3.51)-- (3.45,-4.42);
		\draw [line width=2.8pt] (3.45,-4.42)-- (5.01,-4.1);
		\draw [line width=2.8pt] (5.01,-4.1)-- (5.94,-2.44);
		\draw [line width=1.6pt,dash pattern=on 2pt off 2pt] (1.28,-0.93)-- (1.2,0.53);
		\draw [line width=1.6pt,dash pattern=on 2pt off 2pt] (5.94,-2.44)-- (6.13,-2.07);
		\draw (1.74,-1.25) node[anchor=north west] {$L$};
		\draw (-0.01,0.23) node[anchor=north west] {};
		\begin{scriptsize}
		\draw [fill=qqqqff] (0.58,1.52) circle (1.5pt);
		\draw [fill=qqqqff] (1.78,4.8) circle (1.5pt);
		\draw [fill=qqqqff] (6.3,3.78) circle (1.5pt);
		\draw [fill=qqqqff] (6.48,-2.1) circle (1.5pt);
		\draw [fill=qqqqff] (0.26,-3.1) circle (1.5pt);
		\draw [fill=qqqqff] (-0.44,-0.54) circle (1.5pt);
		\draw[color=black] (6.37,-2.17) node {};
		\end{scriptsize}
		\end{tikzpicture}
		
	\end{center}
	\caption{A sphere with 3 boundary components where the boundary components are in blue. The loop $L$, thick black one, is not a geodesic since it is not length minimizing.}
	\label{convex}
\end{figure}

\begin{remark}
	\label{pants}
	Let $S^L$ be an FDL surface so that $\mathfrak{p},\mathfrak{p'}$ are empty. Assume that it has a boundary component of non-negative curvature.
	Any geodesic loop in $S^L$ either lies in this boundary component or 
	does not intersect with this component. See Figure \ref{convex}. If the loop lies in this boundary component, then this component is non-singular; each point at this  has zero curvature.
\end{remark}

\begin{figure}
	\begin{center}
		\includegraphics[scale=0.5]{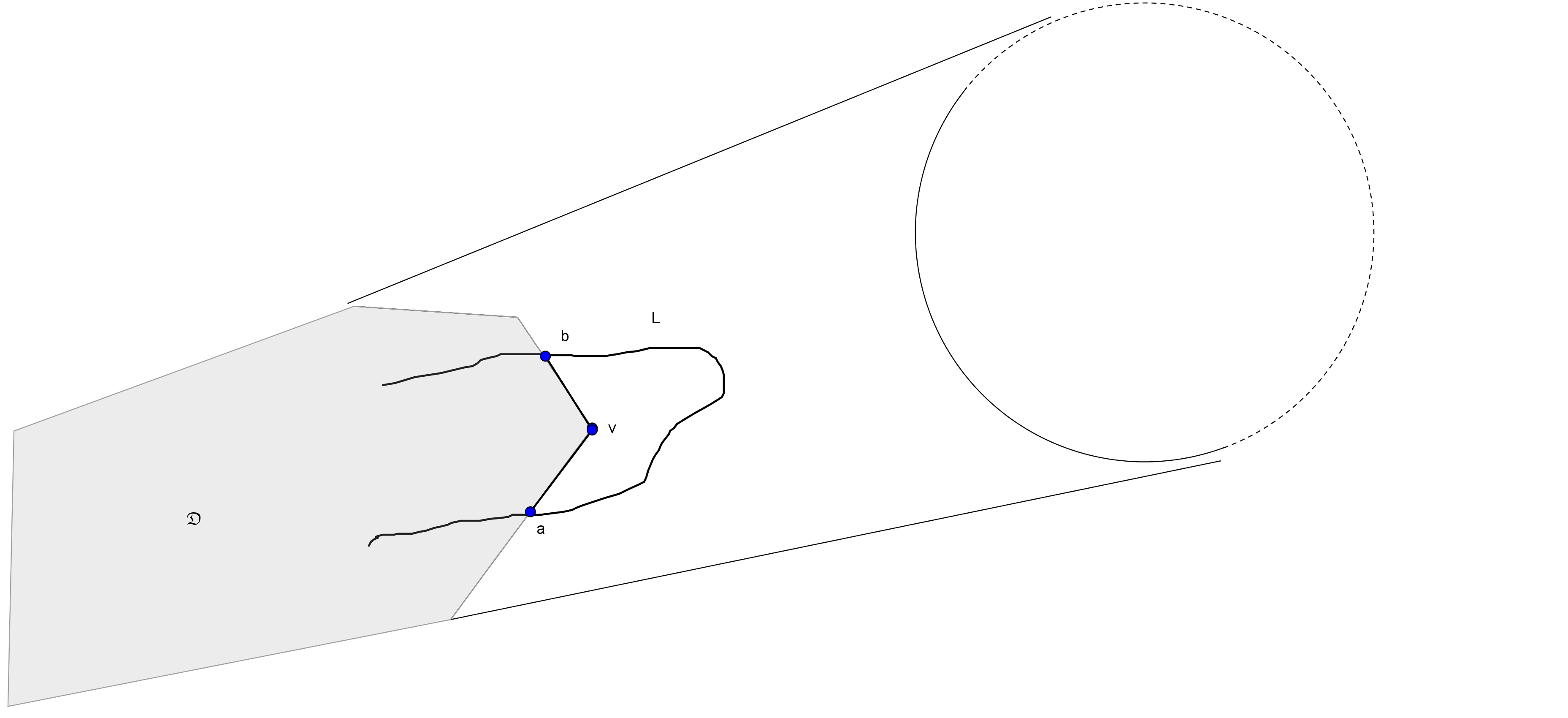}
	\end{center}
	\caption{Since boundary points of $\frak{D}$ have non-negative curvature, a loop which intersects with $D$ and its complement can not be length minimizing.}
	\label{full}
\end{figure}

\begin{remark}
	\label{mainbody}
	Let $S^L$ be a FDL surface. As in the proof of  Proposition \ref{mod1} and Corollary \ref{mod-corollary}, cut $S^L$ through disks around its punctures so that  at each point of each \textit{resulting} boundary component, curvature is non-negative. Let $\frak{D}$ be the resulting compact surface. For each loop $L$ in $S^L$ which intersects such a component, there exists a loop in its homotopy class
	which has length less than or equal to length of $L$ and lies in $\mathfrak{D}$. See
	Figure \ref{full}. The part of the loop $L$ which does not lie in $\frak{D}$ has length greater than $\lvert[a,v]\rvert + \lvert [v,b] \rvert$. This happens since boundary points of $\frak{D}$ has non-negative curvature.

\end{remark}

\begin{theorem}
	\label{freehomo}
	Given an FDL surface $S^L$ and a loop $L$ on it, there exists a length minimizing geodesic on its free homotopy class.
	
	\begin{proof}
		We cut the surface around its punctures as in Lemma \ref{mod1} and Corollary \ref{mod-corollary}  so that resulting disks do not intersect $L$. Thus the part left is 
		a compact surface $\mathfrak{D}$ containing $L$, and  there exists a length minimizing geodesic $g$ in the homotopy class of $L$ in $\mathfrak{D}$. Since the curvature at each point of a \textit{resulting}  boundary component is non-negative  we see that either $g$ lies in such a boundary component  and this boundary component is non-singular, or it does not intersect 
		such a boundary component. See Remark \ref{pants}. This shows that $g$ is indeed a geodesic in $S^L$ and it is in homotopy class of $L$.
		
		Assume that there is a geodesic $g'$ in homotopy class of $g$ whose length is less than length of $g$. It follows that $g'$ does not lie completely in $\frak{D}$, and 
Remark \ref{mainbody} implies that there is a shorter loop $g''$  in same  homotopy class which lies in $\frak{D}$. This implies  that length length of $g''$ is greater than or equal to the length of $g$,
which is a contradiction.  
	\end{proof}
\end{theorem}

\section*{Acknowledgements}
I am really grateful to Muhammed Uluda\u{g} for his suggestions. Also, I thank Deniz Kutluay
for his reading the first version of present 
manuscript. I am grateful to Ayberk Zeytin and 
Susumu Tanabe for their comments.

\label{references}


\end{document}